\documentclass[preprint, prX]{revtex4}

\usepackage{amsmath}    
\usepackage{graphicx}   
\usepackage{verbatim}   
\usepackage{color}      
\usepackage{subfigure}  
\usepackage{hyperref}   
\usepackage{amssymb}    
\usepackage{epsfig}
\usepackage{graphics,graphicx}
\usepackage{setspace}
\usepackage{url}
\usepackage{algorithm,algorithmic}

\begin{document}

\begin{abstract}
Firefly algorithms belong to modern meta-heuristic algorithms inspired by nature that can be successfully applied to continuous optimization problems. In this paper, we have been applied the firefly algorithm, hybridized with local search heuristic, to combinatorial optimization problems, where we use graph 3-coloring problems as test benchmarks. The results of the proposed memetic firefly algorithm (MFFA) were compared with the results of the Hybrid Evolutionary Algorithm (HEA), Tabucol, and the evolutionary algorithm with SAW method (EA-SAW) by coloring the suite of medium-scaled random graphs (graphs with 500 vertices) generated using the Culberson random graph generator. The results of firefly algorithm were very promising and showed a potential that this algorithm could successfully be applied in near future to the other combinatorial optimization problems as well.

\textit{To cite paper as follows: I. Fister Jr., X.-S. Yang, I. Fister, J. Brest, Memetic firefly algorithm for combinatorial optimization,
in Bioinspired Optimization Methods and their Applications (BIOMA 2012), B. Filipic and J.Silc, Eds.
Jozef Stefan Institute, Ljubljana, Slovenia, 2012}

\end{abstract}

\title{Memetic firefly algorithm for\\\ combinatorial optimization}

\author{Iztok Fister Jr.}
\altaffiliation{University of Maribor, Faculty of electrical engineering and computer science
Smetanova 17, 2000 Maribor}
\email{iztok.fister@guest.arnes.si}

\author{Xin-She Yang}
\altaffiliation{University of Cambridge, Department of Engineering 
Trumpington Street, Cambridge CB2 1PZ, UK}
\email{xy227@cam.ac.uk}

\author{Iztok Fister}
\altaffiliation{University of Maribor, Faculty of electrical engineering and computer science
Smetanova 17, 2000 Maribor}
\email{iztok.fister@uni-mb.si}

\author{Janez Brest}
\altaffiliation{University of Maribor, Faculty of electrical engineering and computer science
Smetanova 17, 2000 Maribor}
\email{janez.brest@uni-mb.si}

\maketitle

\section{Introduction}
Nature, especially biological systems, has always been an inspiration for those scientists who would like to transform some successful features of a biological system into computer algorithms for efficient problem solving. Birds, insects, ants, and fish may display some so-called, collective or swarm intelligence in foraging, defending, path finding, etc. This collective intelligence of self-organizing systems or agents has served as a basis for many good and efficient algorithms developed in the past, e.g.: ant-colony optimization \cite{Dorigo:2004}, particle swarm optimization\cite{Kennedy:1995}, artificial bee colony \cite{Fister:2012,Fister:2012a,Karaboga:2007}, bacterial foraging \cite{Liu:2002}. Today, all these algorithms are referred to as \textit{swarm intelligence}.

The firefly algorithm (FFA) belongs to swarm intelligence as well. It was developed by X. S. Yang~\cite{Yang:2008}. This algorithm is based on the behavior of fireflies. Each firefly flashes its lights with some brightness. This light attracts other fireflies within the neighborhood. On the other hand, this attractiveness depends on the distance between the two fireflies. The closer the two fireflies are, the more attractive they will seem to other. In FFA, each firefly represents a point in a search space. When the attractiveness is proportional to the objective function the search space is explored by moving the fireflies towards more attractive neighbors.

FFA has displayed promising results when applied to continuous optimization problems~\cite{Yang:2009,Yang:2010}. Conversely, within the area of combinatorial optimization problems, only a few papers have been published to date. Therefore, aim of this paper is to show that FFA can be applied to this kind of optimization problems as well. In this context, FFA for graph 3-coloring (3-GCP) has been developed. 3-GCP can informally be defined as follows: How to color a graph $G$ with three colors so that none of the vertices connected with an edge is colored with the same color. The problem is \textit{NP}-complete as proved by Garey and Johnson~\cite{Garey:1979}.

The most natural way to solve this problem is in a greedy fashion. Here, the vertices are ordered into a permutation and colored sequentially one after the other. However, the quality of coloring depends on the order in which the vertices will be colored. For example, the \textit{naive} method orders the vertices of graph randomly. One of the best traditional heuristics for graph coloring today is DSatur by Brelaz~\cite{Brelaz:1979}, which orders the vertices $v$ according to \textit{saturation degree} $\rho(v)$. The saturation degree denotes the number of distinctly colored vertices to the vertex $v$.

This problem cannot be solved by an exact algorithm for graph instances of more than 100 vertices. Therefore, many heuristic methods have been developed for larger instances of graphs. These methods can be divided into local search~\cite{Galinier:2006} and hybrid algorithms~\cite{Malaguti:2009}. One of the more successful local search heuristic was Tabucol developed by Herz and De Werra~\cite{Hertz:1987}, who employed the tabu search proposed by Glover~\cite{Glover:1986}. The most effective local search algorithms today are based on reactive partial local search~\cite{Blochliger:2008,Malaguti:2008}, adaptive memory~\cite{Hertz:2008}, and variable search space~\cite{Blochliger:2008}. On the other hand, various evolutionary algorithms have been hybridized using these local search heuristics. Let us refer to three such algorithms only: the hybrid genetic algorithm by Fleurent and Ferland~\cite{Fleurent:1996}, the hybrid evolutionary algorithm by Galinier and Hao~\cite{Galinier:1999}, and the memetic algorithm for graph coloring by L{\"u} and Hao~\cite{Lu:2010}.

Some modifications of the original algorithm need to be performed in order to apply FFA to 3-GCP. The original FFA algorithm operates on real-valued vectors. On the other hand, the most traditional heuristics act on the permutation of vertices. In order to incorporate the benefits of both, solutions of the proposed memetic FFA (MFFA) algorithm are represented as real-valued vectors. The elements of these vectors represent \textit{weigths} that determine how hard the vertices are to color. The higher the weight is, the sooner the vertex should be colored. The permutation of vertices is obtained by sorting the vertices according to their weights. The DSatur traditional heuristic is used for construction of 3-coloring from this permutation. A similar approach was used in the evolutionary algorithm with the SAW method (EA-SAW) of Eiben et al.~\cite{Eiben:1998}, and by the hybrid self-adaptive differential evolution and hybrid artificial bee colony algorithm of Fister et al.~\cite{Fister:2011, Fister:2012}. Additionally, the \textit{heuristical swap} local search is incorporated into the proposed MFFA. In order to preserve the current best solution in the population, the elitism is considered by this algorithm.

The results of the proposed MFFA algorithm for 3-GCP were compared with the results obtained with EA-SAW, Tabucol, and HEA by solving an extensive set of random medium-scale graphs generated by the Culberson graph generator~\cite{Culberson:2012}. The comparison between these algorithms shows that the results of the proposed MFFA algorithm are comparable, if not better, than the results of other algorithms used in the experiments.

The structure of this paper is as follows: In Section 2, the 3-GCP is discussed, in detail. The MFFA is described in Section 3, while the experiments and results are presented in Section 4. The paper is concluded with a discussion about the quality of the results, and the directions for further work are outlined.

\section{Graph 3-coloring}

Let us suppose, an undirect graph $G=(V,E)$ is given, where $V$ is a set of vertices $v \in V$ for $i=1, \ldots, n$, and $E$ denotes a set of edges that associate each edge $e \in E$ for $i=1, \ldots, m$ to the unordered pair $e=\{v_{i},v_{j}\}$ for $i=1, \ldots, n$ and $j=1, \ldots, n$. Then, the vertex 3-coloring (3-GCP) is defined as a mapping $c:V \rightarrow S$, where $S=\{1,2,3\}$ is a set of three colors and $c$ a function that assigns one of the three colors to each vertex of $G$. A coloring $\mathbf{s}$ is \textit{proper} if each of the two vertices connected with an edge are colored with a different color.

3-GCP can be formally defined as a constraint satisfaction problem (CSP). It is represented as a pair $\langle S, \phi \rangle$, where $S$ denotes the search space, in which all solutions $\mathbf{s} \in S$ are \textit{feasible}, and $\phi$ a Boolean function on $S$ (also a \textit{feasibility condition}) that divides the search space into \textit{feasible} and \textit{unfeasible} regions. To each $e \in E$ the constraint $b_{e}$ is assigned with $b_e(\langle s_1, \ldots , s_n \rangle) = true$ if and only if $e=\{v_i,v_j\}$ and $s_i \neq s_j$. Suppose that $B^i=\{b_e|e=\{v_i,v_j\} \wedge j=1 \ldots m\}$ defines the set of constraints belonging to variable $v_i$. Then, the feasibility condition $\phi$ is  expressed as a conjunction of all the constraints $\phi (\mathbf{s})= \wedge_{v \in V} B^{v}(\mathbf{s})$.

As in evolutionary computation, constraints can be handled indirectly in the sense of a \textit{penalty function}, that punishes the unfeasible solutions. The farther the unfeasible solution is from the feasible region, the higher is the penalty function. The penalty function is expressed as:

\begin{equation}
\label{eq:penalty}
 f(\mathbf{s})=min\sum_{i=0}^{n} \psi(\mathbf{s},B^{i}),
\end{equation}

\noindent where the function $\psi(\mathbf{s},B^{i})$ is defined as:

\begin{equation}
\label{eq:viol}
\psi(\mathbf{s},B^{i})=\left\{\begin{matrix}
1 & \textup{if}\ \mathbf{s}\ \textup{violates\ at\ least\ one\ } b_{e} \in B^{i}, \\
0 & \textup{otherwise}.
\end{matrix}\right.
\end{equation}

Note that Eq.~(\ref{eq:penalty}) also represents the objective function. On the other hand, the same equation can be used as a feasibility condition in the sense that $\phi(\mathbf{s})=true$ if and only if $f(\mathbf{s})=0$. If this condition is satisfied a proper graph 3-coloring is found.

\section{Memetic firefly algorithm for graph 3-coloring}

The phenomenon of fireflies is that fireflies flash their lights that can be admired on clear summer nights. This light is produced by a complicated set of chemical reactions. Firefly flashes in order to attract mating partners and serve as a protection mechanism for warning off potential predators. Their light intensity $I$ decreases when the distance $r$ from the light source increases according to term $I \varpropto r^2$. On the other hand, air absorbs the light as the distance from the source increases.

When the flashing light is proportional to the objective function of the problem being optimized (i.e., $I(\mathbf{w}) \varpropto f(\mathbf{w})$), where $\mathbf{w}$ represents the candidate solution) this behavior of fireflies can represent the base for an optimization algorithm. However, artificial fireflies obey the following rules: all fireflies are unisex, their attractiveness is proportional to their brightness, and the brightness of a firefly is affected or determined by the landscape of the objective function.

These rules represent the basis on which the firefly algorithm acts~\cite{Yang:2008}. The FFA algorithm is population-based, where each solution denotes a point in the search space. The proposed MFFA algorithm is hybridized with a local search heuristic. In this algorithm, the solution is represented as a real-valued vector $\mathbf{w_{i}}=(w_{i,1},\ldots,w_{i,n})$ for $i=1 \ldots \mathit{NP}$, where $\mathit{NP}$ denotes the size of population $P$. The vector $\mathbf{\mathbf{w_{i}}}$ determines the weights assigned to the corresponding vertices. The values of the weights are taken from the interval $w_{i,j} \in [lb,ub]$, where \textit{lb} and \textit{ub} are the lower and upper bounds, respectively. The weights represent an initial permutation of vertices $\pi(\mathbf{v_{i})}$. This permutation serves as an input to the DSatur heuristic that obtains the graph 3-coloring. The pseudo-code of the MFFA algorithm is illustrated in Algorithm~\ref{alg:prog}.

\begin{algorithm}[htb]
\caption{Pseudo code of the MFFA algorithm}
\label{alg:prog}
\begin{algorithmic}[1]
\STATE $t$ = 0; $fe$ = 0; $found$ = \textbf{FALSE}; $\textbf{s}^{*} = \emptyset;$
\STATE $P^{(t)}$ = InitializeFFA(); 
\WHILE {(!TerminateFFA($fe$, $found$))}
\STATE $fe$ += EvaluateFFA($P^{(t)}$);
\STATE $P'$ = OrderFFA($P^{(t)}$);
\STATE $found$ = FindTheBestFFA($P^{(t)}, \textbf{s}^{*}$);
\STATE $P^{t+1}$ = MoveFFA($P^{t}, P'$);
\STATE t = t+1;
\ENDWHILE
\end{algorithmic}
\end{algorithm}

As can be seen from Algorithm~\ref{alg:prog}, the search process of the MFFA algorithm (statements within the \textbf{while} loop) that is controlled by generation counter $t$ consists of the following functions:
\begin{itemize}
  \item EvaluateFFA(): evaluating the solution. This evaluation is divided into two parts: In the first part, the solution $\mathbf{w_{i}}$ is transformed into a permutation of vertices $\pi(\mathbf{v_{i})}$ according to the non-decreasing values of the weights. In the second part, the permutation of vertices $\pi(\mathbf{v_{i})}$ is decoded into a 3-coloring $\mathbf{s}_{i}$ by the DSatur heuristic. Note that the 3-coloring $\mathbf{s}_{i}$ represents the solution of 3-GCP in its original problem space, where the quality of the solution is evaluated according to Eq.~(\ref{eq:penalty}). Conversely,  looking for a new solution is performed within the real-valued search space.
  \item OrderFFA(): forming an intermediate population $P'$ by copying the solutions from the original population $P^{(t)}$ and sorting $P^{(t)}$ according to the non-decreasing values of the objective function.
  \item FindTheBestFFA(): determining the best solution in the population $P^{(t)}$. If the best solution in $P^{(t)}$ is worse than the $\textbf{s}^{*}$ the later replaces the best solution in $P^{(t)}$, otherwise the former becomes the best solution found so far $\textbf{s}^{*}$.
  \item MoveFFA(): moving the fireflies towards the search space according to the attractiveness of their neighbour's solutions.
\end{itemize}

Two features need to be developed before this search process can take place: the initialization (function InitializeFFA()) and termination (function TerminateFFA()). The population is initialized randomly according to the following equation:

\begin{equation}
\label{eq:eq1}
 w_{i,j}=(ub-lb) \cdot rand(0,1)+lb,
\end{equation}

\noindent where the function $\mathrm{rand}(0,1)$ denotes the random value from the interval $[0,1]$. The process is terminated when the first of the following two condition is satisfied: the number of objective function evaluations $\mathit{fe}$ reaches the maximum number of objective function evaluations (MAX\_FES) or the proper coloring is found ($\mathit{found} == \mathit{true}$).

The movement of $i$-th firefly is attracted to another more attractive firefly $j$, and expressed as follows:

\begin{equation}
\label{eq:move}
\mathbf{w_{i}}=\mathbf{w_{i}}+\beta_{0} e^{- \gamma r_{i,j}^{2}} (\mathbf{w_{j}}-\mathbf{w_{i}})+ \alpha (rand(0,1)- \frac{1}{2}),\ \mathrm{for}\ j=1...n.
 \end{equation}

\noindent Note that the move of the $i$-th firefly is influenced by all the $j$-th fireflies for which $I[j]>I[i]$. As can be seen from Eq.~\ref{eq:move}, two summands are added to the current firefly position $\mathbf{w_{i}}$. The former reflects the attractiveness between firefly $i$ and $j$ (determined by $\beta(r)=\beta_{0}e^{- \gamma r_{i,j}^{2}}$), while the latter is the randomized move in the search space (determined by the randomized parameter $\alpha$). Furthermore, the attractiveness depends on the $\beta_{0}$ that is the attractiveness at $r=0$, absorbtion coefficient $\gamma$, and Euclidian distance between the attracted and attracting firefly $r_{i,j}$.

\subsection{The heuristical swap local search}

\begin{figure}[htb!]
\vspace{-2mm}
\centering
\includegraphics [scale=0.4]{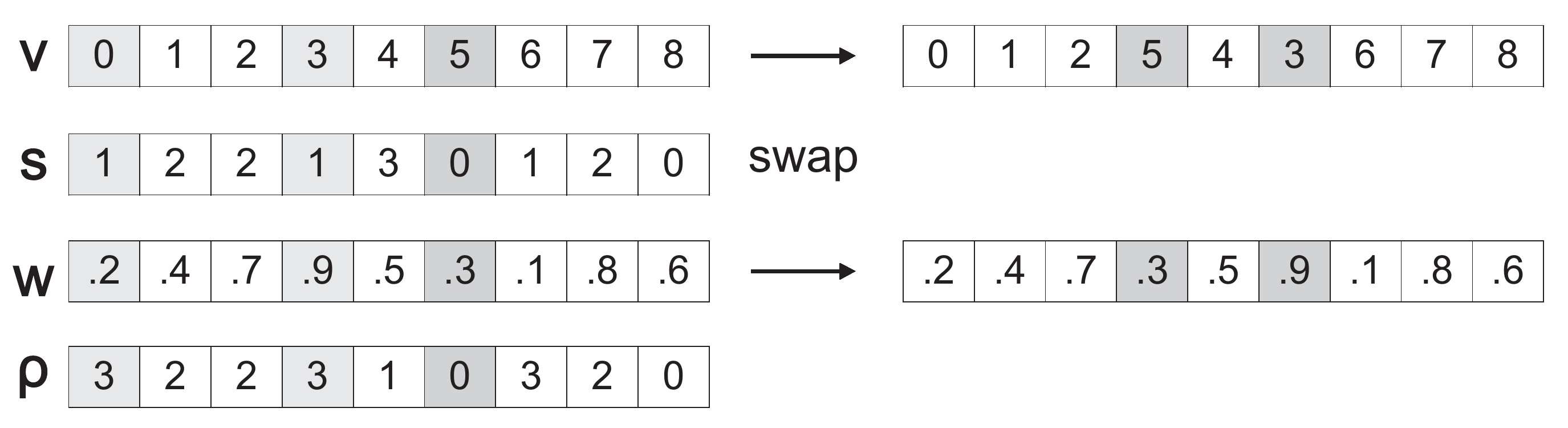}
\caption{The heuristical swap unary operator}
\label{fig:Swap}
\vspace{-2mm}
\end{figure}

After the evaluation step of the MFFA algorithm, the heuristical swap local search tries to improve the current solution. This heuristic is executed until the improvements are detected. The operation of this operator is illustrated in Fig.~\ref{fig:Swap}, which deals with a solution on $G$ with 10 vertices. This solution is composed of a permutation of vertices $\mathbf{v}$, 3-coloring $\mathbf{s}$, weights $\mathbf{x}$, and saturation degrees $\rho$. The heuristical swap unary operator takes the first uncolored vertex in a solution and orders the predecessors according to the saturation degree descending. The uncolored vertex is swapped with the vertex that has the highest saturation degree. In the case of a tie, the operator randomly selects a vertex among the vertices with higher saturation degrees ($1$-opt neighborhood).

In Fig.~\ref{fig:Swap}, an element of the solution corresponding to the first uncolored vertex 5 is in dark gray and the vertices 0 and 3 with the highest saturation degree are in light gray. From vertices 0 and 3, heuristical swap randomly selects vertex 3 and swaps it with vertex 5 (the right-hand side of Fig.~\ref{fig:Swap}).

\section{Results}

In the experimental work, the results of the proposed MFFA algorithm were compared with the results of: EA-SAW, HEA, and Tabucol. The algorithms used in the experiments were not selected coincidentally. In order to help the developers of new graph coloring algorithms, the authors Chiarandini and St{\"u}tzle~\cite{Chiarandini:2010} made the code of Tabucol and HEA available within an online compendium~\cite{Chiarandini:2012}. On the other hand, this study is based on the paper of Eiben and al.~\cite{Eiben:1998}, in which the authors proposed the evolutionary algorithm with SAW method for 3-GCP. The source code of this algorithm can also be obtained from Internet~\cite{Hemert:2012}. The goal of this experimental work was see whether the MFFA could also be applied to combinatorial optimization problems like 3-GCP.

The characteristics of the MFFA in the experiments were as follows. The population size was set at 500. The MAX\_FES was fixed at 300,000 by all algorithms to make this comparison as fair as possible. All algorithms executed each graph instance 25 times. The algorithm parameters of MFFA were set as follows: $\alpha = 0.1$, $\beta_{0} = 0.1$, and $\gamma = 0.8$. Note that these values of algorithm parameters optimize the performance of MFFA and were obtained during parameter tuning within the extensive experimental work. This parameter tuning satisfies the first perspective of parameter tuning, as proposed by Eiben and Smit~\cite{Eiben:2011}, i.e., choosing parameter values that optimize the algorithm's performance. In order to satisfy the second perspective of the parameter tuning, i.e., how the MFFA performance depends on its parameter values, only the influence of the edge density was examined because of limited paper length.

The algorithms were compared according to the measures \textit{success rate} ($SR$) and \textit{average evaluations of objective function to solution} ($AES$). While the first measure represents the ratio between the number of successfully runs and all runs, the second determines the efficiency of a particular algorithm. The aim of these preliminary experiments was to show that MFFA could be applied for 3-GCP. Therefore, any comparison of algorithms according to the time complexity was omitted.

\subsection{Test-suite}

The test-suite, considered in the experiments, consisted of graphs generated with the Culberson random graph generator~\cite{Culberson:2012}. The graphs generated by this generator are distinguished according to type, number of vertices $n$, edge probability $p$, and seeds of random number generator $q$. Three types of graphs were employed in the experiments: \textit{uniform} (random graphs without variability in sizes of color sets), \textit{equi-partite}, and \textit{flat}. The edge probabilities were varied in the interval $p \in 0.008 \ldots 0.028$ with a step of 0.001. Finally, the seeds were varied in interval $q \in 1 \ldots 10$ with a step of one. As a result, $3 \times 21 \times 10 = 630$ different graphs were obtained. That is, each algorithm was executed $15,750$ times to complete this experimental setup.

The experimental setup was selected so that a \textit{phase transition} was captured. The phase transition is a phenomenon that accompanies almost all NP-hard problems and determines the region where the NP-hard problem passes over the state of "solvability" to a state of "unsolvability" and vice versa~\cite{Turner:1988}. Typically, this region is characterized by ascertain problem parameter. This parameter is edge probability for 3-GCP. Many authors have determined this region differently. For example, Petford and Welsh~\cite{Petford:1989} stated that this phenomenon occurs when $2pn/3 \approx 16/3$, Cheeseman et al.~\cite{Cheeseman:1991} when $2m/n \approx 5.4$, and Eiben et al.~\cite{Eiben:1998} when $7/n \leq p \leq 8/n$. In the presented case, the phase transition needed to be by $p=0.016$ over Petford and Welsh, by $p \approx 0.016$ over Cheeseman, and between $0.014 \leq p \leq 0.016$ over Eiben et al.

\subsection{Influence of the edge density}

\begin{figure}[!hbt]	
\vspace{-5mm}
\centering
\subfigure[SR by uniform graphs] {\includegraphics[width=5.9cm]{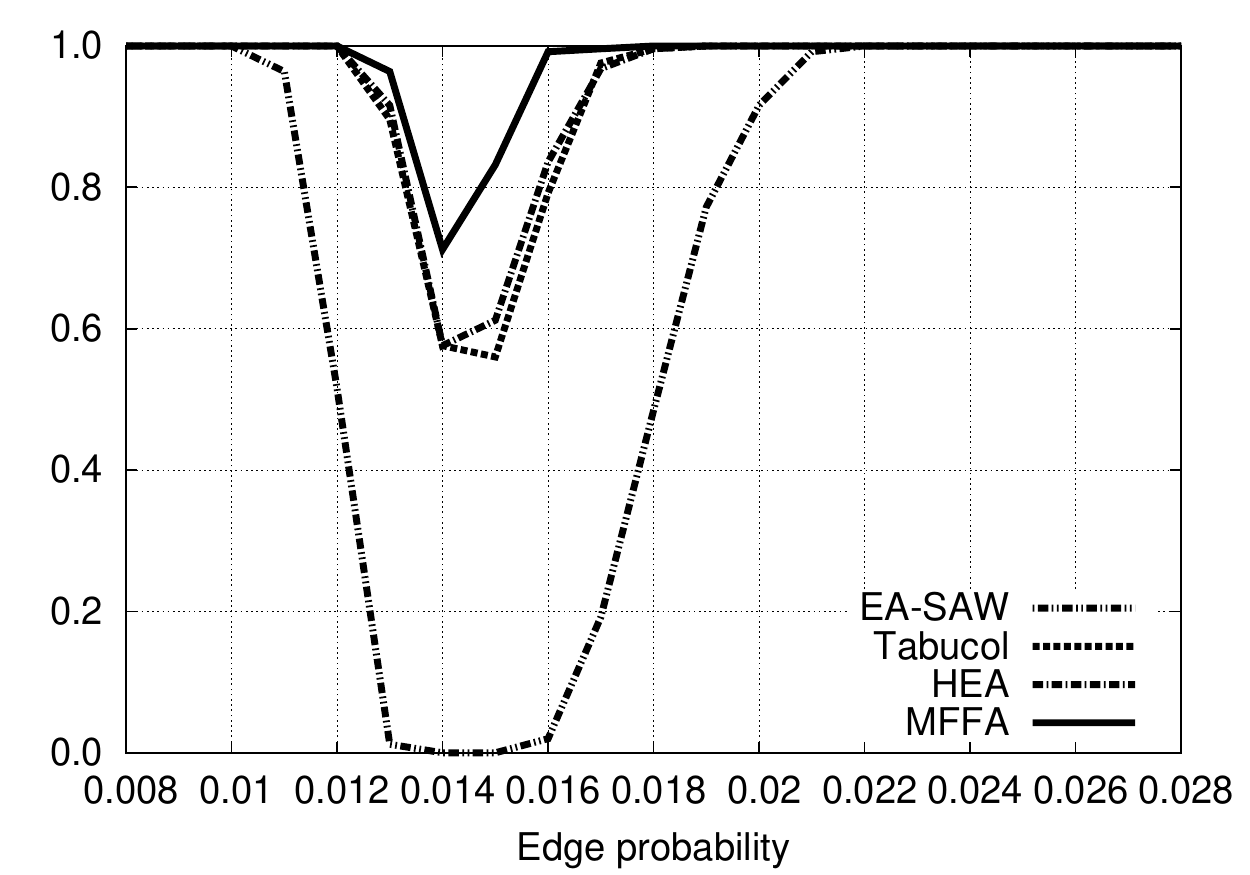}}
\subfigure[AES by uniform graphs] {\includegraphics[width=5.9cm]{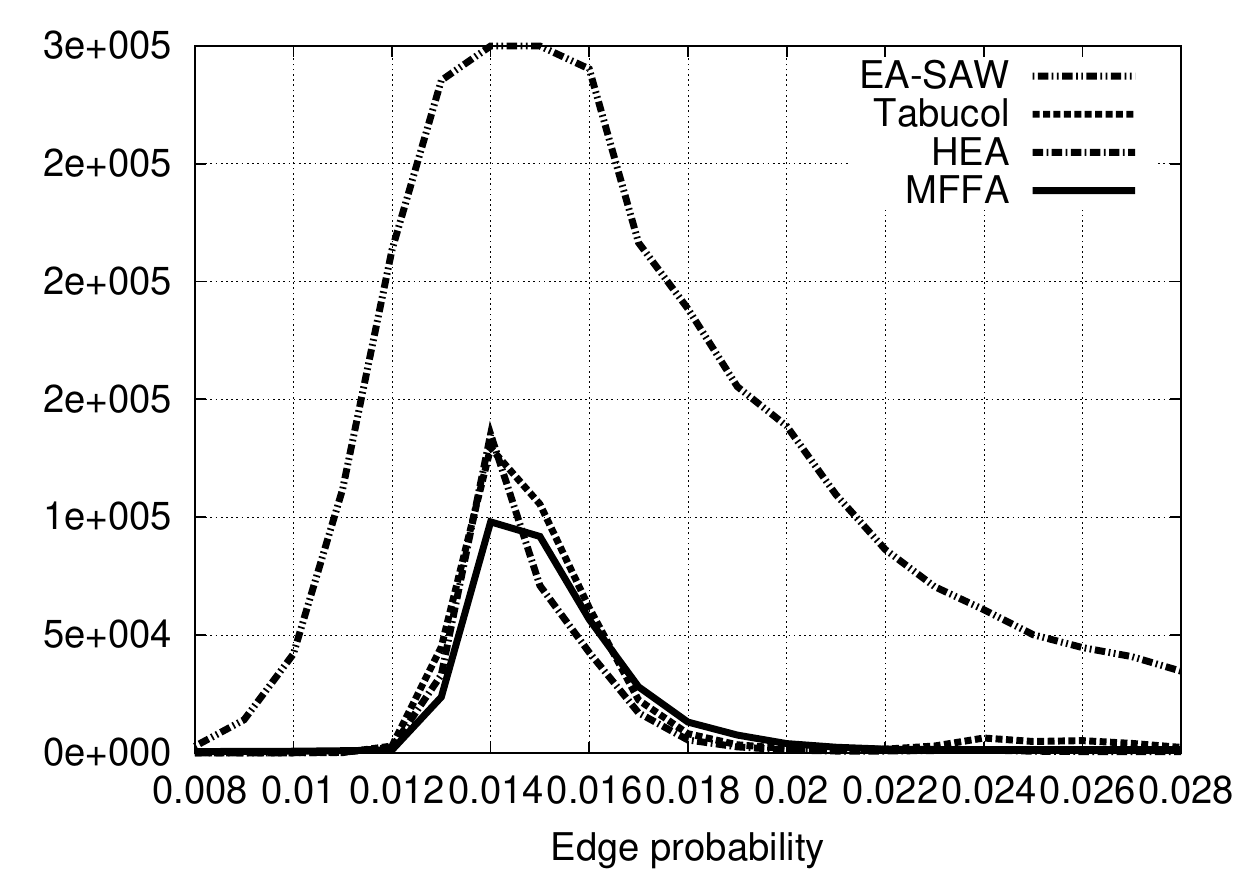}}
\subfigure[SR by equi-partite graphs] {\includegraphics[width=5.9cm]{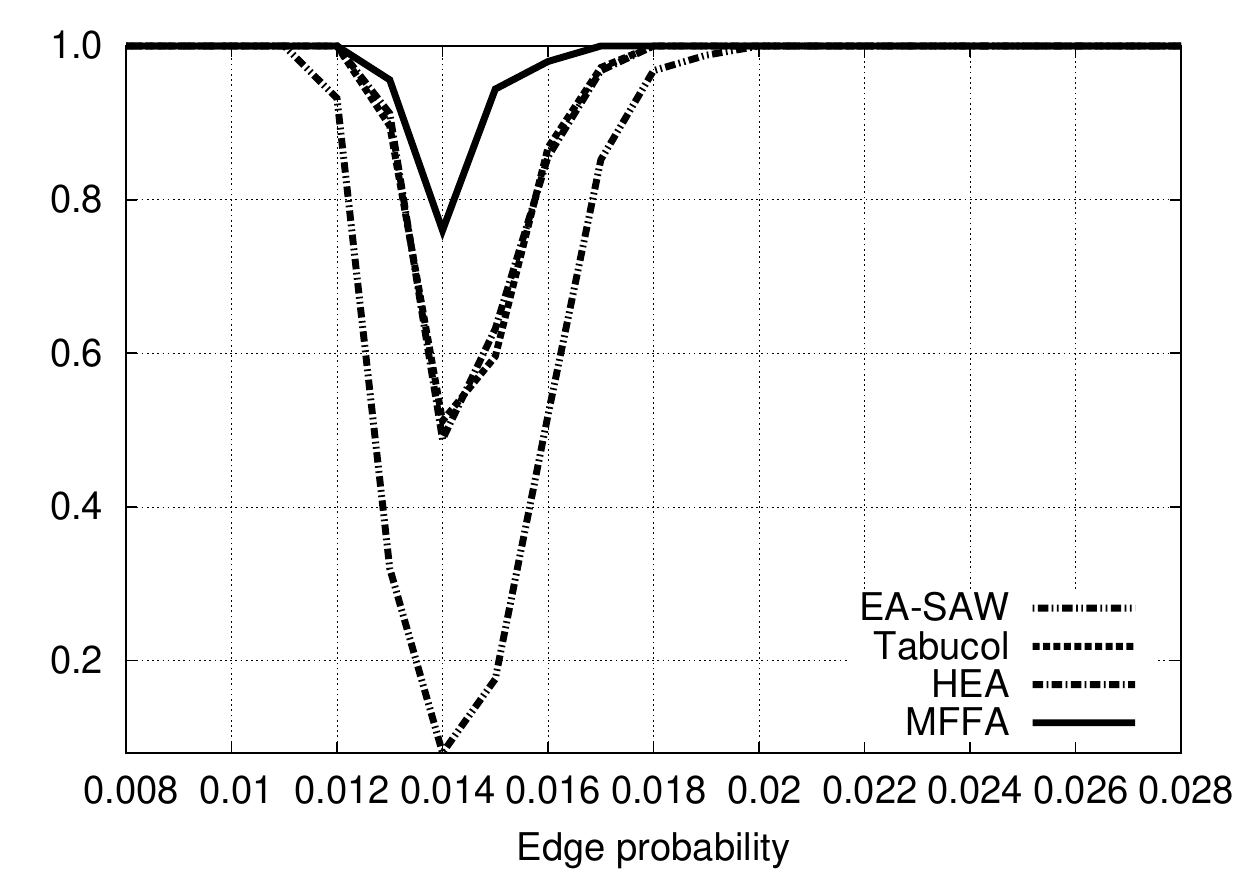}}
\subfigure[AES by equi-partite graphs] {\includegraphics[width=5.9cm]{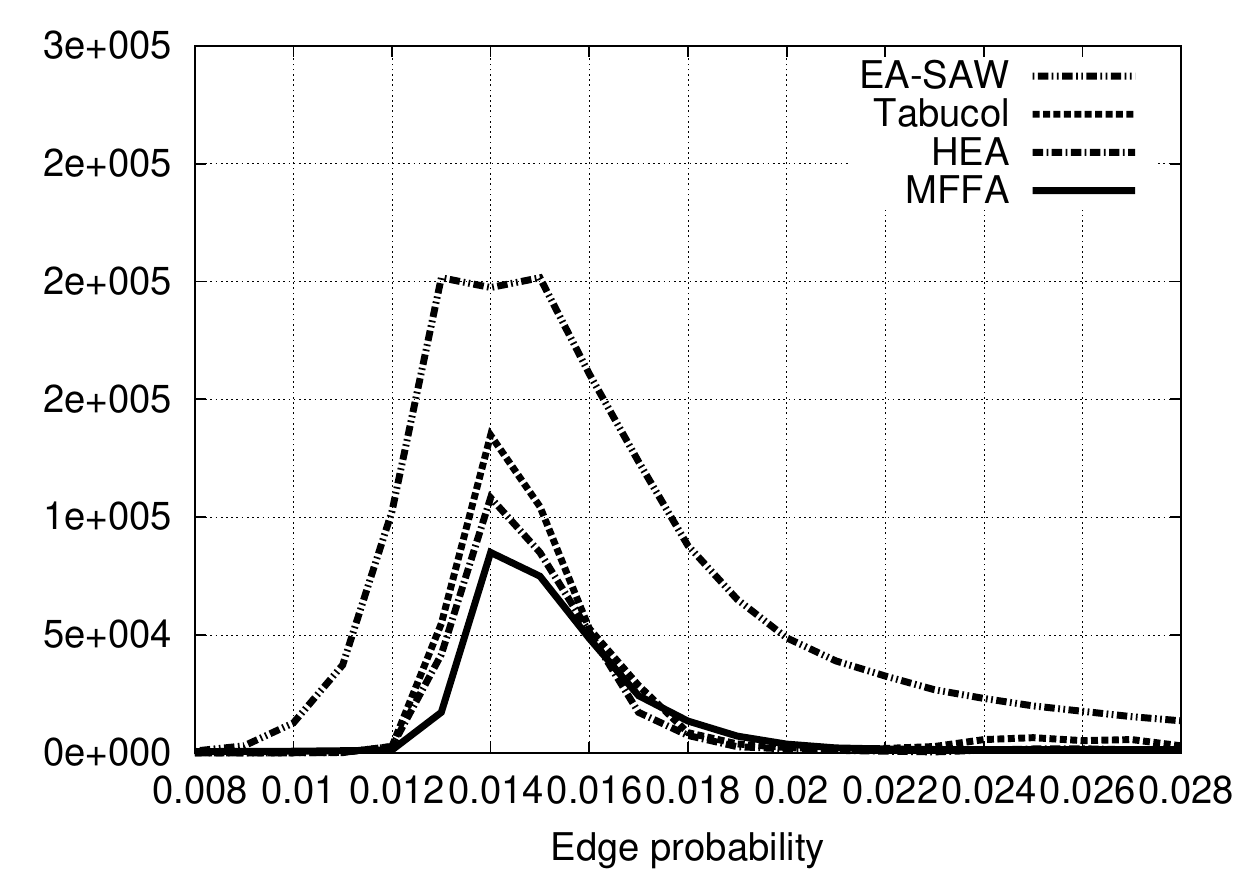}}
\subfigure[SR by flat graphs] {\includegraphics[width=5.9cm]{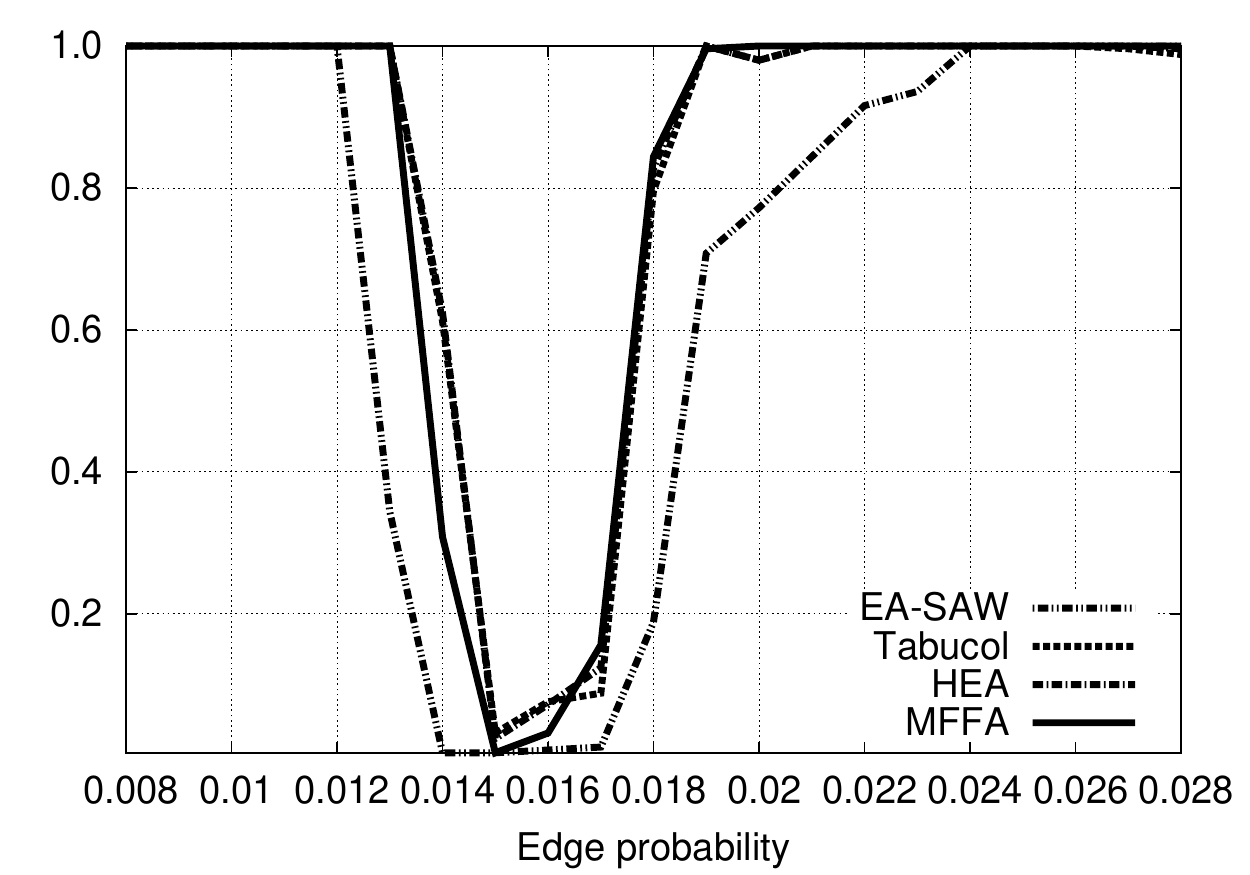}}
\subfigure[AES by flat graphs] {\includegraphics[width=5.9cm]{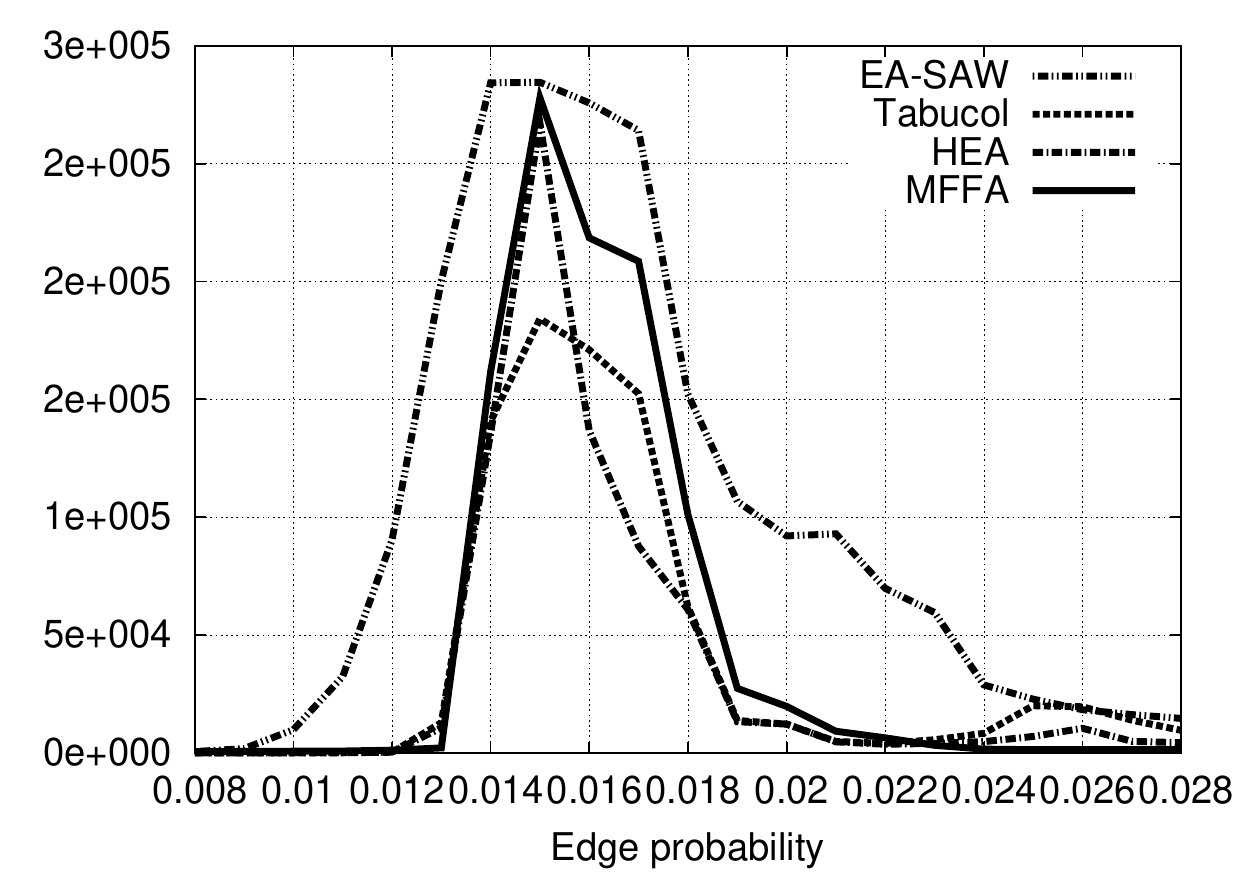}}
\caption{Results of algorithms for 3-GCP solving different types of random graphs}
\label{fig:Sub_1}
\vspace{-5mm}
\end{figure}

During this experiment, the influence of the edge density on the performance of the tested algorithms were investigated. The results are illustrated in Fig.~\ref{fig:Sub_1}. This figure consists of six diagrams corresponding to graphs of different types (uniform, equi-partite, and flat), and according to the measures $SR$ and $AES$. In these diagrams, the average of those values accumulated after 25 runs are presented. Especially, we focus on the behavior of tested algorithms within the phase transition.

As can be seen in Fig.~\ref{fig:Sub_1}.a and Fig.~\ref{fig:Sub_1}.c, the results of MFFA outperformed the results of the other algorithms according to the measure SR on uniform and equi-partite graphs. HEA was slightly better than Tabocol during the phase transition ($p \in [0.014,0.016]$), whilst EA-SAW exposed the worst results within this region. As can be seen in Fig.~\ref{fig:Sub_1}.e, flat graphs were the hardest nuts to crack for all algorithms. Here, the results of MFFA according to measure SR were slightly worse but comparable to the results of HEA and Tabucol, whilst EA-SAW reported the worst results.

According to the measure AES (Fig.~\ref{fig:Sub_1}.b and Fig.~\ref{fig:Sub_1}.d), the best results were reported for MFFA by coloring the uniform and equi-partite graphs. On average, MFFA found solutions using a minimum number of evaluations. Note that the highest peak by $p=0.014$ denotes the hardest instance of uniform and equi-partite graphs to color for almost all the observed algorithms. Conversely, the graph with $p=0.015$ was the hardest instance when coloring the flat graphs.

In summary, the proposed MFFA outperformed the results of HEA and Tabucol when coloring the uniform and equi-partite graphs, while by coloring the flat graphs it behaved slightly worse. The results of EA-SAW fell behind the results of the other tested algorithms for coloring all other types of graphs.

\subsection{Discussion}

The results of other parameter tuning experiments have been omitted because of limited paper length. Notwithstanding, almost four characteristics of MFFA could be exposed from the results of the last experiment. First, it is very important whether the movement of $i$-th firefly according to Eq.~(\ref{eq:move}) is calculated from the position of the $j$-th firefly taken from the population $P^{(t)}$ or the intermediate population $P'$, because the former incorporates an additional randomness within the search process. Thus, the results were significantly improved. Second, an exploration of the search space depends on the best solution in the population that directs the search process to more promising regions of the search space. As a result, this elitist solution needs to be preserved. Third, the local search serves as a search mechanism for detailed exploration of the basins of attraction. Consequently, those solutions that would normally be ignored by the regular FFA search process can be discovered. Fourth, the $\alpha$ parameter determines the size of the randomness move within the search space. Unfortunately, all conducted tests to self-adapt this parameter did not bear any improvement, at this moment. In summary, these preliminary results of MFFA encourage us to continue developing this algorithm in the future.

\section{Conclusion}

The results of MFFA showed that this algorithm could in future be a very promising tool for solving 3-GCP and, consequently, the other combinatorial optimization problems as well. In fact, it produced better results than HEA and Tabucol, when coloring the medium-scale uniform and equi-partite graphs. Unfortunately, this algorithm is slightly worse on flat graphs that remains the hardest to color for all the tested algorithms.

However, these good results could be misleading until further experiments on large-scale graphs (graphs with 1,000 vertices) are performed. Fortunately, in the sense of preserving the obtained results, we have several ways for improving this MFFA in the future.

\bigskip{\small \smallskip\noindent Updated 10 May 2012.}
\end{document}